\documentclass{article}
\usepackage{epsfig}

\title{The Combinatorics of Iterated Loop Spaces.}

\author{M.A. Batanin\protect \footnote{The author holds the Scott Russell Johnson Fellowship in
the Centre of Australian Category Theory at Macquarie University}\\ Macquarie University,  North
Ryde, NSW 2109, Australia \\ e-mail: mbatanin@math.mq.edu.au} 
\date{ }

\newtheorem{theorem}{\bf Theorem}[section]
\newtheorem{defin}{\bf Definition}[section]
\newtheorem{pro}{\bf Proposition}[section]

\newtheorem{cor}{\bf Corollary}[theorem]

\newcommand{\Oper}{\mbox{$ Oper_n^{S}$}}

\newcommand{\h}{\mbox{$\bf h$}}
\newcommand{\Ht}{\mbox{$\bf H$}}
\renewcommand{\H}{
\mbox{{\unitlength=0.25mm
\begin{picture}(13,10)(0,-2.7)
\put(5,2){\makebox(0,0){\mbox{$\cal H$}}}
\put(5.3,2){\makebox(0,0){\mbox{$\cal H$}}}
\put(5.1,2){\makebox(0,0){\mbox{$\cal H$}}}
\put(5,2){\makebox(0,0){\mbox{$\cal H$}}}
\put(4.9,2){\makebox(0,0){\mbox{$\cal H$}}}
\end{picture}}}
}

\renewcommand{\thefigure}{\Roman{figure}}

\newcommand{\V}{
\mbox{$op(V)$}
}

\newcommand{\Top}{\mbox{$op(Top)$}
}

\newcommand{\M}{\mbox{$\cal M$}}

\newcommand{\Proof}{\noindent {\bf Proof. \ }}

\newcommand{\adj}{ {\unitlength=0.25mm
\begin{picture}(30,10)(-10,1)
\put(5,2){\makebox(0,0){\mbox{$\longrightarrow$}}}
\put(5,8){\makebox(0,0){\mbox{$\longleftarrow$}}}
\end{picture}}}

\newcommand{\Q}{
{\unitlength=0.25mm
\begin{picture}(500,10)(-10,0)
\put(440,10){\line(0,-1){10}}
\put(440,0){\line(1,0){10}}
\put(450,0){\line(0,1){10}}
\put(450,10){\line(-1,0){10}}
\put(451,11){\line(-1,0){10}}
\put(451,1){\line(0,1){10}}
\put(450.5,0.5){\line(0,1){10}}
\put(450.5,10.5){\line(-1,0){10}}
\end{picture}}}

\newcommand{\MT}{ {\unitlength=0.71mm
\begin{picture}(25,10)(4,0.5)
\put(15,0){\line(-1,1){5}}
\put(15,0){\line(1,1){5}}
\put(10,5){\line(0,1){5}}
\put(20,5){\line(0,1){5}}

\end{picture}}}   

\newcommand{\MV}{ {\unitlength=0.71mm
\begin{picture}(25,10)(4,0.5)
\put(15,5){\line(-1,1){5}}
\put(15,5){\line(1,1){5}}
\put(15,0){\line(0,1){5}}

\end{picture}}}

\newcommand{\MS}{ {\unitlength=0.71mm
\begin{picture}(25,10)(4,0)
\put(15,0){\line(-1,1){5}}
\put(15,0){\line(1,1){5}}
\put(10,5){\line(1,1){5}}
\put(20,5){\line(0,1){5}}
\put(10,5){\line(-1,1){5}}
\end{picture}}}   

\newcommand{\Example}{\noindent \makebox[24mm]{{\bf Example
\hspace{-1mm}\addtocounter{example}{1} 
\thesection.\theexample \ }}}

\newcommand{\Remark}{\noindent \makebox[23mm]{{\bf Remark
\hspace{-1mm}\addtocounter{remark}{1} 
\thesection.\theremark \ }}}

\newcounter{remark}[section]
\newcounter{example}[section]

\begin{document}

\maketitle

\begin{abstract}

 It is well known since Stasheff's work that 1-fold loop spaces can be described in terms of 
the existence of higher homotopies for associativity (coherence conditions)
 or equivalently  as algebras of  contractible 
non-symmetric operads. The combinatorics of these higher homotopies is well understood 
and is extremely useful. 

 For  $n \ge 2$ the   theory of symmetric operads encapsulated the corresponding higher homotopies,
yet hid the combinatorics and it has  remain a mystery for almost 40 years. However, the recent developments in many fields
ranging from
  algebraic
topology and  algebraic geometry to mathematical physics and category theory show that this combinatorics in higher
dimensions will be  even more important than the one dimensional case. 
 
 In this paper we are going to show that there exists a conceptual way to
 make these combinatorics explicit  using the so called higher 
nonsymmetric $n$-operads.

\end{abstract}

\section{Introduction.}

\subsection{Preoperadic history of the subject.}

For the decade beginning around 1955 the  question of 
characterising of loop spaces through algebraic structures was a hot
 subject in topology.  A very nice solution for 1-fold
loop spaces was provided by J.Stasheff \cite{St1,St2}. Let us sketch it briefly.

Let $X$ be a pointed space with the based point $e$. Suppose also we have a multiplication 
$$\mu:X\times X \rightarrow X,$$
$$\mu(a,b) = ab$$
For simplicity we assume also that $e$ is a two sided unit of this multiplication.
Then the first condition will be the existence of a homotopy
$$\mu(\mu\times 1)\sim \mu(1\times \mu).$$ 
 In other words we have a map
$$\mu_3:K_3\times X^3 \rightarrow X,$$ 
where $K_3$ is the interval $[0,1]$,
which provides a path 
$$a(bc)\  -\!\!\!-\!\!\!-\!\!\!- \ (ab)c$$
for any $a,b,c\in X$.

We should not stop here. If we consider all possible bracketings of $4$ symbols we thus have
a pentagon of homotopies

{\unitlength=1mm

\begin{picture}(65,43)(-25,-2)

\put(30,30){\makebox(0,0){\mbox{$(ab)(cd)$}}}

\put(10,18){\makebox(0,0){\mbox{$((ab)c)d$}}}

\put(50,18){\makebox(0,0){\mbox{$a(b(cd))$}}}

\put(19.5,2.5){\makebox(0,0){\mbox{$(a(bc))d$}}}

\put(40,2.5){\makebox(0,0){\mbox{$a((bc)d)$}}}

\put(15,20){\line(1,1){9}}
\put(36,29){\line(1,-1){9}}
\put(15,16){\line(1,-3){4}}
\put(41,4){\line(1,3){4}}
\put(26,2.2){\line(1,0){7}}

\end{picture}}
\begin{figure}[h]\caption{}\label{assoc}\end{figure} 

Our next condition is: we should be able to extend these homotopies to the map
$$\mu_4: K_4\times X^4 \rightarrow X,$$ 
 where $K_4$ is the pentagon above. It is clear how to proceed now. In general there exists a
sequence of  convex polytopes
$K_n$ for all $n$ ($K_0 = K_1 = K_2 = \star$) called associahedra. The vertices of $K_n$
correspond to all binary bracketings of a string of $n$ letters.

\begin{defin}[Stasheff \cite{St1,St2}] A connected pointed space $X$ with a multiplication $\mu$ is called an
$A_{\infty}$-space
 if  there exists a sequence of continuous maps
$$\mu_n : K_n\times X^n \rightarrow X $$
where $\mu_n$ is an extension of  a map from the boundary of $K_n\times X^n$ which can be
 constructed from low dimensional 
$\mu_d$. \end{defin}
\begin{theorem}[Stasheff \cite{St1,St2}] A connected topological space $X$ is a $1$-fold loop space if and only if
it admits a structure of an $A_{\infty}$-space.
\end{theorem}

Stasheff's approach to recognition of loop spaces turned out to be exceptionally fruitful and was
and still is a source of inspiration for many breakthrough discoveries.

\

 What about double loop spaces? It is clear that the multiplication in a double 
loop space should be homotopy commutative. For  simplicity, let us suppose that
 it is strictly associative. 
 So we have a homotopy:
$$ab \  -\!\!\!-\!\!\!-\!\!\!- \ ba$$
Next, we should have a homotopy filling in the triangle

{\unitlength=0.8mm

\begin{picture}(100,40)(-11,2)

\put(30,32){\makebox(0,0){\mbox{$abc$}}}

\put(15,18){\makebox(0,0){\mbox{$acb$}}}

\put(30,5){\makebox(0,0){\mbox{$cab$}}}

\put(17,20){\line(1,1){9}}

\put(17,16){\line(1,-1){9}}
\put(30,6.5){\line(0,1){23}}


{\unitlength=0.8mm

\begin{picture}(30,40)(-50,0)

\put(30,32){\makebox(0,0){\mbox{$abc$}}}

\put(45,18){\makebox(0,0){\mbox{$bac$}}}

\put(30,5){\makebox(0,0){\mbox{$bca$}}}

\put(43,20){\line(-1,1){9}}

\put(43,16){\line(-1,-1){9}}
\put(30,6.5){\line(0,1){23}}

\put(5,18){\makebox(0,0){\mbox{and also}}}

\end{picture}}
\end{picture}}
\begin{figure}[h]\caption{}\label{triangle}\end{figure}

These homotopies  and the homotopy for commutativity
allow  us to fill in the hexagon in two essentially different ways (compare
with two different methods of proving the Yang-Baxter equation in a
braided monoidal category).

{\unitlength=0.8mm

\begin{picture}(100,50)(-10,-9)

\put(30,32){\makebox(0,0){\mbox{$abc$}}}

\put(20,22.5){\makebox(0,0){\mbox{$acb$}}}

\put(20,10){\makebox(0,0){\mbox{$cab$}}}
\put(30,1){\makebox(0,0){\mbox{$cba$}}}

\put(22,24){\line(1,1){6}}

\put(22,8){\line(1,-1){6}}
\put(20,11){\line(0,1){9.5}}

\put(40,22.5){\makebox(0,0){\mbox{$bac$}}}

\put(40,10){\makebox(0,0){\mbox{$bca$}}}
\put(30,1){\makebox(0,0){\mbox{$cba$}}}

\put(38,24){\line(-1,1){6}}

\put(38,8){\line(-1,-1){6}}
\put(40,11){\line(0,1){9.5}}

\put(29.2,30){\line(-1,-3){6.4}}
\put(37.5,20.5){\line(-1,-3){6.2}}

{\unitlength=0.8mm

\begin{picture}(100,40)(-50,0)

\put(30,32){\makebox(0,0){\mbox{$abc$}}}

\put(20,22.5){\makebox(0,0){\mbox{$acb$}}}

\put(20,10){\makebox(0,0){\mbox{$cab$}}}
\put(30,1){\makebox(0,0){\mbox{$cba$}}}

\put(22,24){\line(1,1){6}}

\put(22,8){\line(1,-1){6}}
\put(20,11){\line(0,1){9.5}}

\put(40,22.5){\makebox(0,0){\mbox{$bac$}}}

\put(40,10){\makebox(0,0){\mbox{$bca$}}}
\put(30,1){\makebox(0,0){\mbox{$cba$}}}

\put(38,24){\line(-1,1){6}}

\put(38,8){\line(-1,-1){6}}
\put(40,11){\line(0,1){9.5}}

\put(30.8,30){\line(1,-3){6.4}}
\put(22.5,20.5){\line(1,-3){6.2}}
\put(5,15){\makebox(0,0){\mbox{and }}}
\end{picture}}

\end{picture}}
\begin{figure}[h]\caption{}\label{YB}\end{figure}

\noindent So we should be able to find a homotopy which fills in the
three dimensional ball with  boundary subdivided according to the
above pictures (S.Crans has considered such an axiom for his theory of
teisi in \cite{CrM}).

{\unitlength=0.8mm

\begin{picture}(100,50)(-35,-9)

\put(30,32){\makebox(0,0){\mbox{$abc$}}}

\put(20,22.5){\makebox(0,0){\mbox{$acb$}}}

\put(20,10){\makebox(0,0){\mbox{$cab$}}}
\put(30,1){\makebox(0,0){\mbox{$cba$}}}

\put(22,24){\line(1,1){6}}

\put(22,8){\line(1,-1){6}}
\put(20,11){\line(0,1){9.5}}

\put(40,22.5){\makebox(0,0){\mbox{$bac$}}}

\put(40,10){\makebox(0,0){\mbox{$bca$}}}
\put(30,1){\makebox(0,0){\mbox{$cba$}}}

\put(38,24){\line(-1,1){6}}

\put(38,8){\line(-1,-1){6}}
\put(40,11){\line(0,1){9.5}}

\put(29.2,30){\line(-1,-3){6.4}}
\put(37,20.5){\line(-1,-3){6.2}}

{\unitlength=0.8mm

\begin{picture}(100,40)(0,0)

\put(30,32){\makebox(0,0){\mbox{$abc$}}}

\put(20,22.5){\makebox(0,0){\mbox{$acb$}}}

\put(20,10){\makebox(0,0){\mbox{$cab$}}}
\put(30,1){\makebox(0,0){\mbox{$cba$}}}

\put(22,24){\line(1,1){6}}

\put(22,8){\line(1,-1){6}}
\put(20,11){\line(0,1){9.5}}

\put(40,22.5){\makebox(0,0){\mbox{$bac$}}}

\put(40,10){\makebox(0,0){\mbox{$bca$}}}
\put(30,1){\makebox(0,0){\mbox{$cba$}}}

\put(38,24){\line(-1,1){6}}

\put(38,8){\line(-1,-1){6}}
\put(40,11){\line(0,1){9.5}}

\multiput(30.8,30)(2.5,-7.5){3}{\line(1,-3){1.5}}
\multiput(22.5,20.5)(2.3,-6.9){3}{\line(1,-3){1.5}}

\end{picture}}

\end{picture}}
\begin{figure}[h]\caption{}\label{YB2}\end{figure}

 If we try to proceed further in this direction  the situation
quickly becomes unmanageable. Stasheff himself noticed that there should be
 a way to write down these conditions explicitly but `it is tediously difficult' \cite{St3}. 

\subsection{ Operads and homotopy coherence.} 

To handle this situation the language of operads was invented by M.Boardman, R.Vogt  and P.May \cite{BV,May}.

\begin{defin} A nonsymmetric (topological) operad is a sequence of topological spaces 
$$A_0, A_1, A_2, \ldots $$
(called the {\it underlying collection}) together with unit $e\in A_1$ and multiplication
$$\mu: A_k\times A_{i_1}\times\ldots\times A_{i_k} \rightarrow A_{i_1+\ldots+i_k}$$
which satisfies associativity and unitarity conditions. 
\end{defin}
The nonsymmetric operads form a category. The morphisms are levelwise morphisms of underlying 
collections which preserve multiplication and units.

There are two important examples. The first is the endomorphism operad $E(X)$ of a topological space
$X$: 
$$E(X)_n = Top(X^n,X).$$
The unit is given by the identity morphism and multiplication by iterated composition of 
$n$-ary maps. The second example is the sequence of Stasheff's 
associahedra. Here multiplication is given by an appropriate inclusion of
 $K_k\times K_{i_1}\times\ldots\times K_{i_k}$ to the boundary of  $K_{i_1+\ldots+i_k}$
defined by Stasheff \cite{St1,St2}. 

\begin{defin} An algebra of an operad $A$ is a topological space $X$ equipped with a morphism of operads
$$A\rightarrow End(X).$$
\end{defin}

The Stasheff's recognition principle for a topological space $X$ can be formulated now as the
existence of a $K$-algebra structure on $X$. Moreover, one can use  any
other contractible nonsymmetric operad for the recognition principle. The operad $K$ is then an
initial (up to higher homotopies) object in the category of contractible nonsymmetric operads (if
we ignore units). So the language of nonsymmetric operads is esentially equivalent to the direct 
combinatorial approach of
Stasheff.

However, for the study of homotopy commutativity the nonsymmetric operads are not enough. 
\begin{defin}
 A symmetric (topological) operad is a sequence of topological spaces 
$$A_0, A_1, A_2, \ldots $$
(called the {\it underlying collection}) together with unit $e\in A_1$ and multiplication
$$\mu: A_k\times A_{i_1}\times\ldots\times A_{i_k} \rightarrow A_{i_1+\ldots+i_k}$$
and actions of the symmetric groups $S_n$ on $A_n$ 
which satisfy associativity, unitarity and equivariancy conditions. 
\end{defin}

 The
endomorphism operad has a natural action of the symmetric groups and can be completed to a
symmetric operad structure. It is now easy to define algebras of a symmetric operads.

\begin{theorem}[May \cite{May}] A connected topological space $X$ is an infinite loop space
if and only if it has an action of a contractible symmetric operad. 
\end{theorem}

 For  finite values of $n$ there are symmetric operads which detect $n$-fold loop spaces: the
so called little $n$-cube operads ${\cal C}^n$ and anything equivalent to them ({\it $E_n$-operads}).  Yet, they are very
complicated from a homotopy point of view: ${\cal C}_k^n$ is homotopy equivalent to the space of
configurations of $k$-distinct points in $n$-dimensional  real space. {\it There is no known
characterisation of $E_n$-operads similar to the characterisation of $E_1$ and $E_{\infty}$
operads except for $n=2$ } \cite{F}.

Now, if we try to write down explicitly the coherence laws for an $n$-fold loop space using the
action of an $E_n$-operad we will have the same trouble as in the previous subsection. 
So the approach based on symmetric operads, being  very powerful in many respects, still
tells almost nothing about the combinatorics of higher homotopies for $n>1$. We quote J.Baez and P.May:
`\ldots work of Boardman and Vogt, May, and Segal gave conceptual encapsulations that hid the 
implicit higher homotopies, whose
combinatorial structure is still somewhat obscure' \cite{BaezMay}.

\subsection{The approach of this paper.}

The main goal of the paper is to give a conceptual approach to  the
combinatorics of $n$-fold loop spaces. It has two essential ingredients.

First we suggest to replace symmetric operads by    higher nonsymmetric operads  which 
appeared in higher category theory \cite{BM}. The intuition behind this approach is that
$E_n$-space
 must be equivalent to {\it $n$-tuple weak $\omega$-groupoid} i.e. a weak
$\omega$-groupoid with one object, one arrow, \ldots , one $(n-1)$-arrow, 
as have been conjectured  by J.Baez
and J.Dolan \cite{BaezD}.
 The decisive step to make this intuition precise was made by us in \cite{BEH}. 
Here  we already  proved   that in this theory we  have a simple characterisation of so called $(n-1)$-terminal
$n$-operads (see Definition \ref{terminal}) which detect $n$-fold loop spaces: they are just the contractible $n$-operads.

The second ingredient is the idea to consider some sort of `universal categorical model' for a given 
theory of operads. To explain this better let us recall a well known fact that the poset of faces of the
associahedron can be described in terms of planar trees. Moreover, the planar trees can be organised into a category,
and even  a categorical  nonsymmetric operad $\bf k$ such that the classifying space of ${\bf k}$ is exactly $K$ \cite{Oper}. The
classical tree formalism is based on the observation that every operad produces a  functor on ${\bf k}$. In \cite{BEH}
we observed that the existence of this formalism is due to the fact that  the  categorical operad  $\bf k$ is actually a
 categorical nonsymmetric operad (without nullary operations) freely generated by an internal nonunital 
nonsymmetric operad. Actually, in \cite{BEH} we work with the full operadic structure with nullary opertions and units and we
construct a categorical symmetric operad $\h^1$ freely generated by an internal 
nonsymmetric operad. It is easy to make appropriate modifications if we want to produce ${\bf k}$. 

In some sense the combinatorics of $E_1$-spaces is governed by the categorical operad ${\bf k}$ or $\h^1$ if we want to 
work with the full structure including homotopy units.  

Using these ideas   we show that
 there exist  appropriate $n$-dimensional analogues of $\h^1$  
and the combinatorics of
$E_n$-spaces is governed by this categorical $n$-operad. Even though this
combinatorics is much more complicated than in Stasheff's case we show that there is no obstacle
to writing it down explicitly.

We also consider the relationships between  $n$-operads and symmetric operads more closely.
In \cite{BEH} we established  an adjunction
 $$Des_n:SOper(W)  \adj Oper_n(\Sigma^n W) : EH_n$$
 between these two categories of operads, which is  the usual adjunction between 
symmetric and nonsymmetric operads in a symmetric monoidal category $W$ when $n=1$. In this paper we  define model category
structures on the categories  of topological $(n-1)$-terminal
$n$-operads and topological symmetric operads  and show that the
adjunction  above gives rise to a
Quillen adjunction \cite{H}. 

Using this techniques we demonstrate that the classifying space $B(\h^n)$ of
the categorical operad $\h^n$ constructed in
\cite{BEH} is  
 a cofibrant replacement for $B(\M^n)$, where $\M^n$ is the operad 
of iterated monoidal categories  introduced in \cite{BFSV}.

This allows us to look closer at the combinatorial structure of $E_n$-operads and we show that
such an operad is actually homotopically freely generated by its internal $n$-operad. 

 In the last section we illustrate on some examples how this
theory can be applied to produce some coherence laws. We are going to  develop this theme in a further paper.

\section{Higher operads.}
\subsection{Trees, their morphisms and $n$-operads.}

In this section we give a brief overview of the theory developed in \cite{BM,BEH} to make our paper
relatively self contained. We do not give  full definitions of monoidal globular category and $n$-operad here
because we   need only some properties of them established in \cite{BM,BEH}. 

 For a natural number $n$ we will denote by $[n]$ the ordinal
$$1 \ < \ 2 \ < \ \ldots \ < n .$$
In particular $[0]$ will denote the empty ordinal. 

\begin{defin} A tree of height $n$ (or simply $n$-tree) is a chain of  order preserving maps of
ordinals
$$T=
[k_n]\stackrel{\rho_{n-1}}{\longrightarrow}[k_{n-1}]\stackrel{\rho_{n-2}}{\longrightarrow}...
\stackrel{\rho_0}{\longrightarrow} [1]
$$  
 \end{defin}
If $i\in [k_m]$ and there is no $j\in [k_{m+1}]$ such that $\rho_{m}(j)=i$ then we call $i$ {\it
a leaf of $T$ of height $i$}.   We will call the leaves of $T$
of  height 
$n$ {\it the tips} of 
$T$. If for an $n$-tree $T$ all its leaves are tips we call
such a tree {\it pruned}.

We illustrate the definition in a picture

{\epsfxsize=250pt 
\makebox(290,110)[r]{\epsfbox{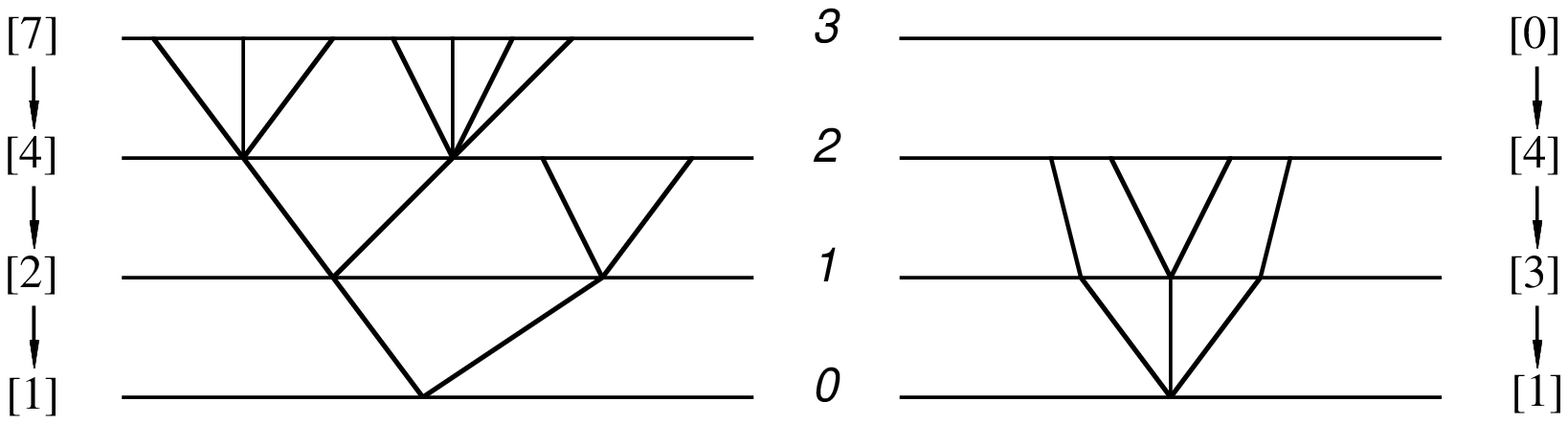}}}
\begin{figure}[h]\caption{}\label{trees}\end{figure}

\noindent The tree on the left side of the picture is not pruned since it has two leaves which are not tips.
The tree on the right side  has the empty ordinal at the highest level; we will
call such trees {\it degenerate}. There is actually an operation on trees which we denote by
$z(-)$ which assigns to the $n$-tree
$[k_n]\stackrel{}{\rightarrow}[k_{n-1}]\stackrel{}{\rightarrow}...
\stackrel{}{\rightarrow} [1]$ the $(n+1)$-tree
$$[0]\longrightarrow
[k_n]\stackrel{}{\longrightarrow}[k_{n-1}]\stackrel{}{\longrightarrow}...
\stackrel{}{\longrightarrow} [1] .$$

Another operation $\partial(-)$ on trees is {\it truncation}
$$\partial([k_n]\stackrel{}{\rightarrow}[k_{n-1}]\stackrel{}{\rightarrow}...
\stackrel{}{\rightarrow} [1])= ([k_{n-1}]\stackrel{}{\rightarrow}...
\stackrel{}{\rightarrow}) [1] .$$

We now define both the source and target of a tree $T$ to be equal to $\partial( T)$. So we have a globular
structure on the set of all trees. We actually have more. The  trees form an $\omega$-category
$Tr$ with the set of $n$-cells being equal to the set of the trees of height $n$. 
If two $n$-trees $S$ and $T$ have the same $k$-sources and $k$-targets (i.e.
$\partial^{n-k}T=\partial^{n-k}S$ ) then they can be composed, and the  composite will be
denoted by $S\otimes_k T$. Then  $z(T)$ is  the operation of taking the identity of
the $n$-cell $T$. Here is an example of the $2$-categorical operations on trees

{\epsfxsize=270pt 
\makebox(300,130)[r]{\epsfbox{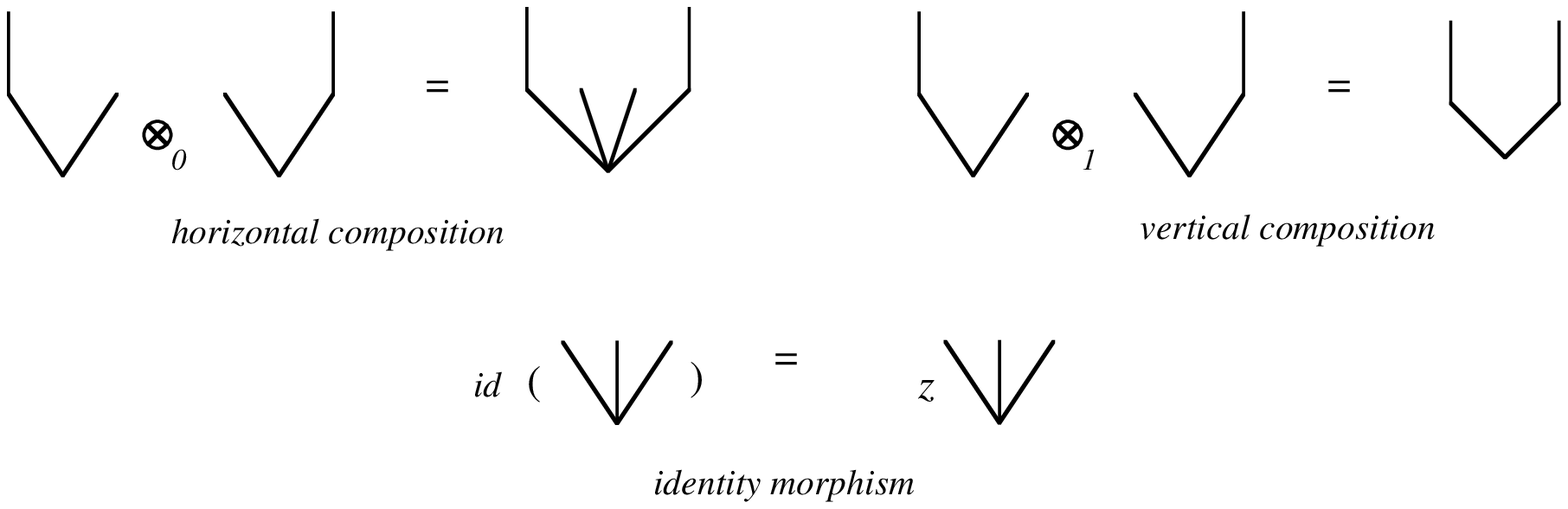}}}
\begin{figure}[h]\caption{}\label{catttrees}\end{figure}

The $\omega$-category $Tr$ is actually the free $\omega$-category generated by the terminal globular
set. Every $n$-tree can be considered as a special sort of $n$-pasting diagram called {\it
globular.} This construction was called the $\star$-construction in \cite{BM}. Here are a couple
of examples.

 {\epsfxsize=200pt 
\makebox(280,150)[r]{\epsfbox{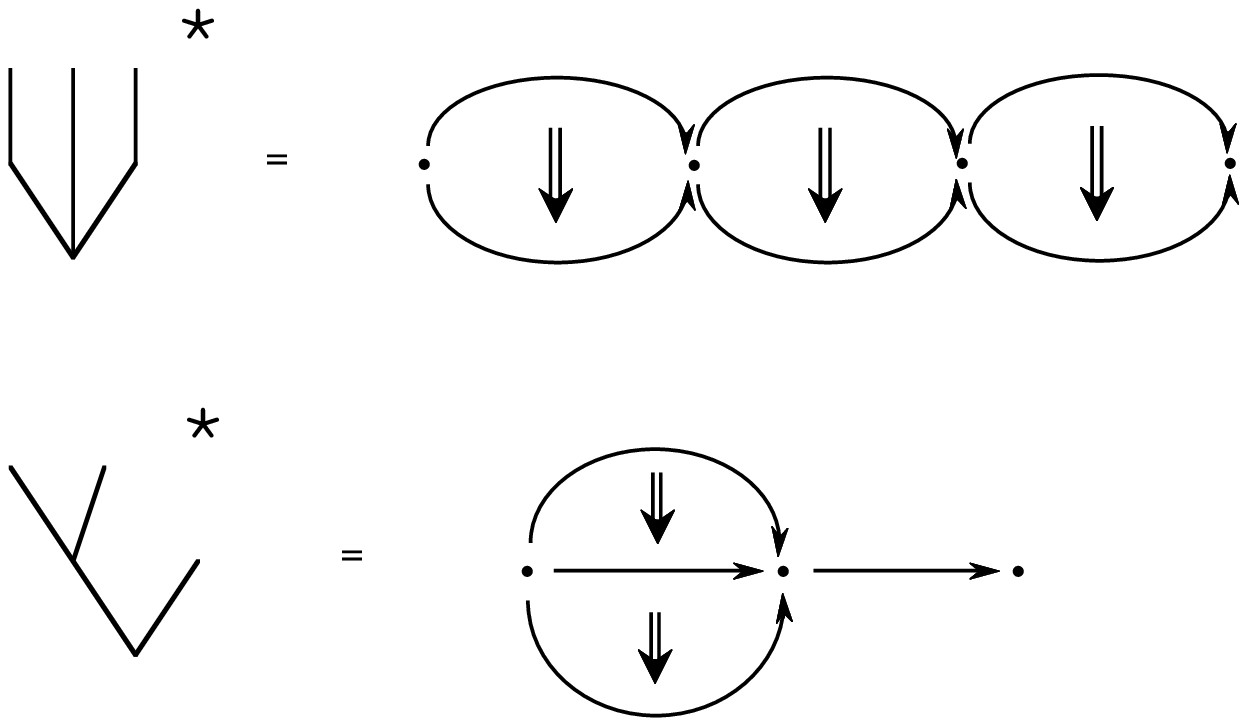}}}
\begin{figure}[h]\caption{}\label{globpast}\end{figure}

 For a globular set $X$ one can then form  the set $D(X)$
of all globular pasting diagrams labelled in $X$. This
is  the free $\omega$-category generated by $X$.
In this way we have a monad $(D,\mu,\epsilon)$ on the
category of globular sets, which plays a central role in
\cite{BM}. 

In particular, $D(1)= Tr.$ We also can consider
$D(Tr) = D^2(1)$. It was observed in \cite{BM} that the $n$-cells of $D(Tr)$ can be identified
with the morphisms of another category introduced by A.Joyal in \cite{J}. This category was called
$\Omega_n$  and some properties of it were studied  later in  \cite{BSt}. More precisely,
it was found that the collection of categories $\Omega_n$ forms an $\omega$-category in $Cat$ and, moreover, it
is freely generated by an internal $\omega$-category. So it is a higher dimensional analogue of the algebraic simplicial
category $\Delta$ (which is of course the free monoidal category generated by a monoid \cite{ML}).  The $\Omega_n$ will
be of primary importance for us.

\begin{defin} The category $\Omega_n$ has as objects the trees of height $n$. The morphisms of $\Omega_n$ are
commutative diagrams    

{\unitlength=0.9mm

\begin{picture}(60,35)(-12,0)

\put(10,25){\makebox(0,0){\mbox{$ [k_n]$}}}
\put(10,20){\vector(0,-1){10}}
\put(12,15){\shortstack{\mbox{$ $}}}

\put(17,25){\vector(1,0){10}}
\put(17,5){\vector(1,0){10}}
\put(18,26){\shortstack{\small\mbox{$\rho_{n-1} $}}}
\put(18,6){\shortstack{\small\mbox{$\xi_{n-1} $}}}

\put(35,25){\makebox(0,0){\mbox{$ [k_{n-1}]$}}}
\put(35,20){\vector(0,-1){10}}
\put(35,5){\makebox(0,0){\mbox{$ [s_{n-1}]$}}}

\put(42,25){\vector(1,0){10}}
\put(42,5){\vector(1,0){10}}
\put(43,26){\shortstack{\small\mbox{$\rho_{n-2} $}}}
\put(43,6){\shortstack{\small\mbox{$\xi_{n-2} $}}}

\put(56,25){\shortstack{\mbox{$.  \   .  \  .  $}}}
\put(56,5){\shortstack{\mbox{$ . \ . \ . $}}}

\put(80,25){\makebox(0,0){\mbox{$[1] $}}}
\put(80,20){\vector(0,-1){10}}

\put(66,25){\vector(1,0){10}}
\put(66,5){\vector(1,0){10}}
\put(69,26){\shortstack{\small\mbox{$\rho_0 $}}}
\put(69,6){\shortstack{\small\mbox{$\xi_0 $}}}

\put(57,21){\shortstack{\mbox{$ $}}}

\put(57,28){\shortstack{\mbox{\small $ $}}}

\put(10,5){\makebox(0,0){\mbox{$ [s_n]$}}}

\put(80,5){\makebox(0,0){\mbox{$[1] $}}}

\put(57,8){\shortstack{\mbox{\small $ $}}}

\put(5,14){\shortstack{\mbox{\small $\sigma_n $}}}
\put(26.5,14){\shortstack{\mbox{\small $\sigma_{n-1} $}}}
\put(75,14){\shortstack{\mbox{\small $\sigma_0 $}}}

\end{picture}}

\noindent where vertical maps are not necessary order preserving but for all $i$ and all $j\in [k_{i-1}]$ the
restriction of
$\sigma_i$ on  $\rho^{-1}_{i-1}(j)$ preserves the natural order on it.
\end{defin}

The category $\Omega_n$ has terminal object 
$$ U_n = [1]\longrightarrow \ldots \longrightarrow [1].$$
Let $T$ be an $n$-tree and let $i$ be a leaf of height $m$ of $T$. Then $i$ determines a unique morphism 
$\xi_i: z^{n-m}U_m\rightarrow T$ in
$\Omega_n$ such that $\xi_{m}(1)=i$. We will often identify the leaf with this morphism. 

Let $\sigma:T\rightarrow S$ be a morphism in $\Omega_n$ and let $i$ be a leaf of $T$ . Then 
{\it the fiber  of $\sigma$ over $i$} is the following pullback in $\Omega_n$ 

{\unitlength=0.9mm

\begin{picture}(60,35)(-2,-1)

\put(45,25){\makebox(0,0){\mbox{$ \sigma^{-1}(i)$}}}
\put(45,22){\vector(0,-1){14}}

\put(77,25){\makebox(0,0){\mbox{$z^{n-m}U_m $}}}
\put(75,22){\vector(0,-1){14}}

\put(54,25){\vector(1,0){13}}

\put(77,14){\shortstack{\mbox{\small $\xi_i $}}}

\put(45,5){\makebox(0,0){\mbox{$ T$}}}

\put(75,5){\makebox(0,0){\mbox{$S $}}}

\put(52,5){\vector(1,0){18}}

\put(60,6){\shortstack{\small \mbox{$\sigma $}}}

\put(57,8){\shortstack{\mbox{\small $ $}}}

\end{picture}}

\noindent which can be calculated as a levelwise pullback in $Set$.

Now, for such a $\sigma:T\rightarrow S$ one can construct a labelling of the
pasting scheme $S^{\star}$  in the $\omega$-category $Tr$ by associating to a
vertex $i$ from $S$ the fiber of $\sigma$ over  $i$. The result of the
pasting operation will be exactly $T$. The inverse process also works i.e. every pasting diagram of trees 
determines a unique morphism in $\Omega_n$. Because of this we can specify a morphism of trees
by the list of its fibers. For example, the two  diagrams

  {\unitlength=0.71mm
\begin{picture}(100,30)(-15,0)
\put(15,5){\line(-1,1){5}}
\put(15,5){\line(1,1){5}}
\put(15,0){\line(0,1){5}}

\put(9,11){\line(1,1){5}}
\put(4,16){\line(0,1){5}}
\put(9,11){\line(-1,1){5}}

\put(21,11){\line(1,1){5}}
\put(26,16){\line(0,1){5}}
\put(21,11){\line(-1,1){5}}

\put(105,5){\line(-1,1){5}}
\put(105,5){\line(1,1){5}}
\put(105,0){\line(0,1){5}}

\put(99,11){\line(1,1){5}}
\put(104,16){\line(0,1){5}}
\put(99,11){\line(-1,1){5}}

\put(111,11){\line(1,1){5}}
\put(106,16){\line(0,1){5}}
\put(111,11){\line(-1,1){5}}

\end{picture}}
\begin{figure}[h]\caption{}\label{twodiagram}\end{figure}

\noindent determine two morphisms in 
$\Omega_2$ from \MT to \MV.
  At level $2$ the first morphism is an identity while the second is the switch map.
Analogously we can present a chain of morphisms 
$T_1\stackrel{\sigma_1}{\rightarrow}T_2 \stackrel{\sigma_2}{\rightarrow}\ldots \stackrel{\sigma_{m-1}}{\rightarrow} T_m$ by
specifying fibers of $\sigma_{m-1}$, fibers of $\sigma_{m-2}$ over fibers of $\sigma_{m-1}$ and so on.

Now an $n$-operad in an augmented monoidal $n$-globular category $V$ can be defined as a collection of
objects $A_{T}, \ T\in Tr_k \ , \ 0\le k\le n $, which respects source and target functors ({\it $n$-collection}),
together with an appropriate unit  and multiplication which make $A$ a monoid in the monoidal category of $n$-collections
\cite{BM}.

In any augmented monoidal $n$-globular category $V$ there exists an operad $\mbox{{\sc T}}(V)$ which algebras are {\it
$n$-globular monoids in}
$V$
\cite{BSt}. In a particular case $V = Span(C)$, where $C$ is cartesian closed category, this operads is just the terminal operad
(and the $n$-globular monoids are internal $n$-categories in $C$).

A funny but important example of an $n$-operad  in the augmented monoidal globular category $Span(CAT)$ ({\it
categorical $n$-operad for simplicity}) is  the following. 
Let $V$ be a strict augmented $n$-globular monoidal category \cite{BM}.
Then we put 
$${\V}_T = V_n, \ T\in Tr_n.$$
The multiplication in $\V$ is given by iterated tensor products in
$V$ and the unit object is given by $z^k (I)$ where $I$ is the unit of
$V_0$. 
  
\subsection{Tree formalism.}

In \cite{BEH} we showed that the tree formalism for operads  is just another way to represent the corresponding category of
operads.

There exists a categorical  $n$-operad $\H^n$ which represents the category of  $n$-operads 
in the following sense. 

There are  two inverse natural  
isomorphisms of categories
$$Oper_n(V) \adj CATOper_n(\H^n, \V) .$$ 
 \noindent Here $Oper_n(V)$ is the category of $n$-operads in $V$ and $CATOper_n(A,B)$ is the category
of operadic functors from a categorical $n$-operad $A$ to a categorical $n$-operad $B$ and their operadic
natural transformations. 

We can  describe the operad $\H^n$ explicitly. An  object of arity $T\in Tr_n$ of $\H^n$ is a 
chain of morphisms
$$\delta = (T\rightarrow T_1 \rightarrow \ldots \rightarrow T_m\rightarrow U_n) $$ 
in $\Omega_n$. There is a morphism $\delta\rightarrow \delta'$ if $\delta'$ can be obtained from $\delta$ by a chain of
composites of  several morphisms in the chain $\delta$ and of insertion  identity morphisms.   
The category $\H^n_{T}$ is contractible since $(T\rightarrow U_n)$ is the terminal
object of it.

Every augmented monoidal globular categry $V$ has its  $(n-1)$-trancation $tr_{n-1}V$ defined by 
$(tr_{n-1}V)_k = V_k, \ 0\ge k\le n-1$. Analogously we can trancate an $n$-operad by considering its restriction
on the trees of height less or equal to $n-1$.

\begin{defin}\label{terminal}
An $n$-operad is called $(n-1)$-terminal if its
$(n-1)$-trancation is $\mbox{{\sc T}}(tr_{n-1} V)$ \cite{BEH}.
\end{defin}

\Remark If $V= \Sigma(W)$ for a braided monoidal category $W$ 
the categories of $1$-operads in $V$ and  of $0$-terminal $1$-operads in $V$ coincide and are isomorphic to the
category of nonsymmetric operads in $W$ \cite{BEH}. If $V = Span(C)$ and  $n=1$ the category of $0$-terminal $1$-operads in
$V$ is the same as the category of nonsymmetric operads in $C$.   However, in general the category of $1$-operads is larger
than the category of $0$-terminal $1$-operads. 

In \cite{BEH} we  introduced the categorical $n$-operad $\Ht^n$
which represents   $(n-1)$-terminal
$n$-operads in the same sense as $\H^n$ represent all  $n$-operads.  The operad $\Ht^n$ is
also a contractible $n$-operad. Let $V$ have globular colimits \cite{BM} and let $ Oper^{(n-1)}_n(V)$ denote the category of
$(n-1)$-terminal
$n$-operads in
$V$. Then we have a pair of adjoint functors
$$\tau_n: Oper^{(n-1)}_n(V) \adj Oper_n(V): \lambda_n,$$
where the right adjoint $\tau_n$ is the obvious inclusion. 

Also in \cite{BEH} we constructed a categorical  symmetric operad $\h^n$ which represents the category of 
$n$-operads in a symmetric categorical operad. The operad $\h^n$ has the homotopy type of the little $n$-cube operad.
We also constructed a pair of adjoints $Des_n \vdash EH_n$:
 $$Des_n:SOper(W)  \adj Oper_n(\Sigma^n W) : EH_n$$
which generalizes the classical adjunction between symmetric and nonsymmetric operads
$$SOper(W) \adj Oper_1(\Sigma W).$$
Here $SOper(W)$ is the category of symmetric operads in a symmetric monoidal category $W$,
 
In the particular case $V= Span(C)$ \cite{BM}, where $C$ is a cartesian closed category, we can identify 
$Oper_n(\Sigma^n C) $ with $Oper^{(n-1)}_n(Span(C))$ so we have a chain of adjunctions
$$SOper(C)\adj Oper^{(n-1)}_n(Span(C)) \adj Oper_n(Span(C)).$$
Moreover, we have a theorem
\begin{theorem}[\cite{BEH}] 
The functors $\tau_n$ and $Des_n$ preserve endomorphism operads, so 
the category of $(n-1)$ terminal algebras of an
$n$-operad
$A$ (i.e. algebras with terminal $(n-1)$-skeleton)
 is isomorphic to the category of algebras of the symmetric operad $EH_n(\lambda_n(A))$. \end{theorem}

Finally, it was proved in \cite{BEH}[Example 10.2] that 
\begin{equation}\label{hn} EH_n(\Ht^n) \simeq \h^n ,\end{equation}
and argument uses only the universal properties of $\Ht^n$ and $\h^n$. Analogously, one can easily show that 
\begin{equation}\label{Hn}
\lambda(\H^n)\simeq \Ht^n .
\end{equation}

\Remark According to (\ref{hn}) the operad $\h^1$ is just a symmetrisation of the nonsymmetric operad $\Ht^1$. 
Yet, the operad $\H^1$ is different from $\Ht^1$. Roughly speaking, $\Ht^1$ is the classical tree operad (if we ignore
units) with morphisms being contractions of internal edges \cite{Oper}.  The objects of the operad $\H^1$ are `regular trees'
in the sense that they are constructed from corollas 
{\epsfxsize=70pt 

\makebox(300,70)[b]{\epsfbox{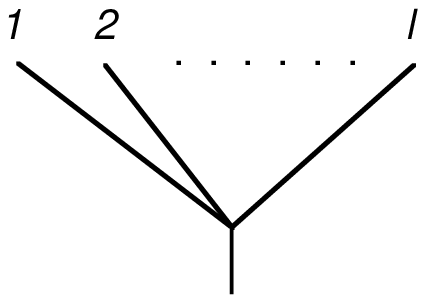}}}
\begin{figure}[h]\caption{}\label{corolla}\end{figure}

\noindent by {\bf simultaneous } grafting at every vertex. The morphisms are
simultaneous contraction of all internal edges on a given level. Every  `regular tree' has its length and so we
have a map
$l$ from the set of regular trees to the set of natural numbers. The monoid of natural numbers is actually  a part of the
structure of $\H^1$; it is  the category 
$\H^1_{U_0}$. The map $l$ is the map of source and target in $\H^1$.  The substitution in $\H^1$ is
 possible only if source is
equal to target and, therefore,  does not destroys regularity of trees.

\section{Model structure on the category of $(n-1)$-terminal $n$-operads}

In this section we adopt the theory of \cite{BergerM} to the case of $(n-1)$-terminal $n$-operads.
 Let $V$ be a cartesian closed model category
 that  is cofibrantly generated and has the terminal object  cofibrant, and 
let $V$ have a symmetric monoidal fibrant
replacement functor. Then the category of $(n-1)$-terminal $n$-collections has an obvious fiberwise 
 model structure. So we define  a fibration
(weak equivalence) of $(n-1)$-terminal
$n$-operads to be an operadic morphism such that its underlying morphism is a 
fibration (weak equivalence) of $n$-collections. 
The argument from \cite{BergerM} works without changes, and we have:
\begin{pro} The category of $(n-1)$-terminal $n$-operads in $V$ is a cofibrantly generated model category.
\end{pro}

Of special interest for us is $V=Top$, the category of compactly generated 
Hausdorf topological. The weak equivalences are weak homotopy equivalences and the fibrations are Serre
fibrations. In this case we can define the bar-construction $B(F,F,X)$ for an $n$-operad $X$, where $F$ is the
free
$(n-1)$-terminal
$n$-operad functor. 

\begin{theorem} 
Let $X$ be an  $(n-1)$-terminal $n$-operad with  cofibrant underlying $n$-collection. Then 
the canonical operad morphism
$$B(F,F,X)\longrightarrow X$$
is  a  cofibrant replacement for $X$.
\end{theorem}

\Proof We have to prove that $B(F,F,X)$ is a cofibrant $n$-operad. 
Let $f:E\rightarrow B$ be a trivial fibration of $n$-operads.

We have to prove that  any operadic map $B(F,F,X)\rightarrow B$ can be lifted to $E$

{\unitlength=1mm

\begin{picture}(60,30)

\put(69,7){\makebox(0,0){\mbox{$B$}}}
\put(69,22){\vector(0,-1){12}}

\put(69,25){\makebox(0,0){\mbox{$E $}}}

\put(55,7){\vector(1,0){10}}

\put(71,15){\shortstack{\mbox{$f $}}}

\put(45,7){\makebox(0,0){\mbox{$B(F,F,X)$}}}

\put(48,10){\line(3,2){5}}
\put(54,14){\line(3,2){5}}
\put(60,18){\vector(3,2){5}}

\end{picture}}

\noindent By construction this amounts to the following lifting problem in the
category of cosimplicial spaces

{\unitlength=1mm

\begin{picture}(60,30)

\put(69,7){\makebox(0,0){\mbox{$\Oper(F^{\star}X,B)$}}}
\put(69,22){\vector(0,-1){12}}

\put(69,25){\makebox(0,0){\mbox{$\Oper(F^{\star}X,E)$}}}

\put(43,7){\vector(1,0){12}}

\put(71,15){\shortstack{\mbox{$f^{\star} $}}}

\put(40,7){\makebox(0,0){\mbox{$\Delta$}}}

\put(43,10){\line(3,2){5}}
\put(49,14){\line(3,2){5}}
\put(55,18){\vector(3,2){5}}

\end{picture}}
  
\noindent 
Here $\Delta$ is the cosimplicial simplicial set consisting of standard simplices and 
 $\Oper$ means the simplicial enriched $Hom$-functor
on the category of $(n-1)$-terminal $n$-operads. 
Since  $\Delta$ is cofibrant in the Reedy model structure \cite{H} it remains to
show that
$f^{\star}$ is a trivial fibration. 

 We follow  a method developed in \cite{BA}. We have
to prove that in the diagram

{\unitlength=1mm

\begin{picture}(100,45)(-10,0)

\put(69,7){\makebox(0,0){\mbox{$M_i(\Oper(F^{\star}X,B))$}}}
\put(69,22){\vector(0,-1){12}}

\put(69,25){\makebox(0,0){\mbox{$M_i(\Oper(F^{\star}X,E))$}}}

\put(40,7){\vector(1,0){12}}

\put(71,15){\shortstack{\mbox{$M_i f^{\star} $}}}

\put(8,35){\vector(1,-2){12.5}}
\put(22,35){\vector(4,-1){30}}
\put(14,35){\vector(1,-1){7}}

\put(25,7){\makebox(0,0){\mbox{$\Oper(F^{i+1}X,B)$}}}

\put(12,38){\makebox(0,0){\mbox{$\Oper(F^{i+1}X,E)$}}}

\put(25,25){\makebox(0,0){\mbox{$W_i$}}}
\put(25,22){\vector(0,-1){12}}
\put(30,25){\vector(1,0){22}}

\put(30.5,20){\line(-1,0){4}}
\put(30.5,20){\line(0,1){3.5}}

\put(22,31){\makebox(0,0){\mbox{$\omega_i$}}}

\end{picture}}

\noindent the canonical map to the pullback is a trivial
fibration. In this diagram $M_i(-)$ is the $i$-th  matching object of
the corresponding cosimplicial object \cite{H}. 

According to Lemma 2.3 from \cite{BA} the diagram above is isomorphic
to the diagram

{\unitlength=1mm

\begin{picture}(100,45)(-10,0)

\put(69,7){\makebox(0,0){\mbox{$Call^S_n(L_{i-1} F^{\star-1}X,B)$}}}
\put(69,22){\vector(0,-1){12}}

\put(69,25){\makebox(0,0){\mbox{$Call^S_n(L_{i-1} F^{\star-1}X,E)$}}}

\put(38,7){\vector(1,0){12}}

\put(71,15){\shortstack{\mbox{$M_i f^{\star} $}}}

\put(8,35){\vector(1,-2){12.5}}
\put(22,35){\vector(4,-1){28}}
\put(14,35){\vector(1,-1){7}}

\put(25,7){\makebox(0,0){\mbox{$Call^S_n(F^{i}X,B)$}}}

\put(12,38){\makebox(0,0){\mbox{$Call^S_n(F^{i}X,E)$}}}

\put(25,25){\makebox(0,0){\mbox{$W_i$}}}
\put(25,22){\vector(0,-1){12}}
\put(28,25){\vector(1,0){22}}

\put(30.5,20){\line(-1,0){4}}
\put(30.5,20){\line(0,1){3.5}}

\put(22,31){\makebox(0,0){\mbox{$\omega_i$}}}

\put(40,33){\makebox(0,0){\mbox{$\phi_i$}}}
\put(12,21){\makebox(0,0){\mbox{$\psi_i$}}}

\end{picture}}

\noindent Here, $Call^S_n$ means the simplicial enriched $Hom$-functor
on the category of $n$-collections, $L_i(-)$ is the latching object for augmented (!) cosimplicial
objects
\cite{H} and $\phi_1, \psi_i$ are generated by the canonical morphism
$$\lambda_{i-1}: L_{i-1} F^{\star-1}X\rightarrow F^i X .$$
If this morphism were a cofibration, then $\omega_i$ would be a trivial fibration by
the axiom for simplicial model category. 

We actually will prove that  $\lambda_{i-1}$ is an isomorphism onto a summand.
For an $n$-collection $X$, let $\tilde{X}$ be an operadic functor from $\Ht^n_d$ to $\Top$
where  $\Ht^n_d$ is a discretisation of the categorical operad $\Ht^n$. Recall \cite{BEH}
that  $\Ht^n_d$ is a free operad generated by the terminal collection, so every topological
$n$-collection considered as a morphism $1\rightarrow \Top$ of $n$-collections  indeed
generates an operadic functor from  $\Ht^n_d$.

Then it is not hard to check that the  free $n$-operad
functor on an $n$-collection can be defined by the following formula:
\begin{equation}\label{freeop}
F(X)_T = \coprod_{T\leftarrow W}\tilde{X}(W)
\end{equation}
where the coproduct is taken over all morphisms in $\Ht^n_T$ (recall that $T$ is the terminal
object of $\Ht^n_T$ ).  It follows that the augmented cosimplicial space
$(F^{\star-1}X)_T$ in dimension $i$ is just the  coproduct
$$\coprod_{T\leftarrow W_0\leftarrow \ldots \leftarrow W_i}\tilde{X}(W_i).$$
 The coface operators are canonical inclusions on the summands corresponding to the
operators of insertion of the identities to the chain 
 $$T\leftarrow W_0\leftarrow \ldots \leftarrow W_i.$$ 
The rest of the proof follows in  complete analogy with Lemma 4.1 of \cite{BCC}.

\Q

Recall \cite{BEH} that on the level of  collections the functor $Des_n$ is defined as
$$Des_n(X)_T = X_{|T|}, $$
where $|T|$ is the number of tips of $T$.
Therefore, $Des_n$ preserves fibrations and weak equivalences and so $Des_n \vdash EH_n$ is a Quillen
adjunction. So  the functor $EH_n$ preserves cofibrations and, in particular, it maps  cofibrant $n$-operads to 
cofibrant symmetric operads.

\begin{cor} The operad $B(\Ht^n)$ is a cofibrant $(n-1)$-terminal $n$-operad. The operad $N(\h^n)$ is a cofibrant 
replacement for the operad $B(\M^n)$ of \cite{BFSV}. 
\end{cor}

\Proof It is not hard to see from the formula (\ref{freeop}) that the operad $B(\Ht^n)$ is isomorphic to
$B(F,F,1)$. But $B(\h^n) \simeq EH_n (B(\Ht^n))$ by (\ref{hn}) and so it is cofibrant. The Theorem 9.2 from \cite{BEH} states
that the the nerve of the canonical morphism $\h^n\rightarrow \M^n $ is a trivial fibration.

\Q

\begin{cor} The restriction of the total left derived functor of $EH_n$
 to the subcategory of contractible $n$-operads induces
an equivalence between the homotopy category of $E_n$-operads and homotopy category of
contractible $(n-1)$-terminal $n$-operads. \end{cor}

\begin{theorem} The homotopy category of $E_n$-spaces  is equivalent to the following three categories
\begin{itemize}
\item 
the homotopy category of
$B(\h^n)$-algebras;
\item 
the homotopy category of
$B(\Ht^n)$-algebras;
\item 
the homotopy category of $(n-1)$-terminal 
$B(\H^n)$-algebras.
\end{itemize}
 \end{theorem}

\section{Internal $n$-operads.}

In this section we show that the categorical theory  of internal operads 
 developed in \cite{BEH} has its topological analogue. Here "$n$-operad" always means "$(n-1)$-terminal
$n$-operad".

\begin{defin} An internal $n$-operad   in a topological symmetric operad $A$ is a coherent
 $n$-operadic functor
$$1 \rightarrow Des_n(A);$$ 
or equivalently, it  is an operadic map $$B(F,F,1)\rightarrow Des_n(A).$$
\end{defin}

For any symmetric topological operad $C$ consider the space ${\mbox{\sc Oper}}_n(C)$ of all internal operads in $C$ 
i.e. the simplicial set 
 $${\mbox{\sc Oper}}_n(C) = \Oper(B(F,F,1),Des_n(C)).$$ 
Observe that the operad $B(\h^n)$ contains a canonical internal $n$-operad given by the unit 
$$B(F,F,1)\rightarrow Des_n(EH_n(B(F,F,1)).$$

\begin{theorem} Let $C$ be an $E_n$-operad. Then $C$ contains an internal $n$-operad. \end{theorem}

\Proof  Since $B(\h^n)$ is a cofibrant $E_n$-operad  we have an operadic equivalence 
$$B(\h^n)\rightarrow C.$$ 
Since $B(\h^n)$ contains an $n$-operad so does $C$.

\Q

\begin{theorem}\label{representable} The functor ${\mbox{\sc Oper}}_n(C)$ is representable on the category of topological
symmetric operads. The representing object is $B(\h^n)$.
\end{theorem}

\Proof Obviously the adjunction isomorphism
$$(\Oper)_0(B(F,F,1),Des_n(C)) \simeq (SOper^S)_0(EH_n(B(F,F,1)),C),$$
where $SOper^S$ is the simplicially enriched $Hom$-functor on the category of symmetric operads,
can be extended to the simplicial adjunction
$$ \Oper(B(F,F,1),Des_n(C))\simeq SOper^S(EH_n(B(F,F,1)),C)\simeq$$
$$\ \ \ \ \ \ \ \ \ \ \ \ \ \ \ \ \ \ \ \ \  \simeq SOper^S(B(\h^n),C)$$

\Q

\Remark It is instructive to try to write down an internal $n$-operad in the little $n$-cube operad ${\cal C}^n$. 

 First let us
choose an orientation of the space
$R^n$.  For a tree
$T$ we can produce a subdivision of the unit
$n$-cube $I^n$ in the following way. If $T= T_1\otimes_0 \ldots \otimes_0 T_k$ is the canonical decomposition
of
$T$ \cite{BEH} then we subdivide $I^n$ into $k$ parallelepipeds  by $k-1$ hyperplanes $$x_1= i/k, \ 1\le i \le
k-1 .$$  Then we consider the canonical decomposition of $T_i = T_{i1}\otimes_1\ldots \otimes_1\ldots T_{il} $
and the corresponding subdivision of the $i$-th parallelepiped by hyperplanes 
$$x_2= i/l, \ 1\le i \le l-1 .$$   
We proceed by induction and get the required subdivision of $I^n$. The interiors of these
parallelepipeds labelled by  natural order on  tips of $T$ gives a configuration of little $n$-cubes $a_T$ which we
take as the vertex of the internal operad corresponding to the tree $T$. Figure (\ref{Milgram}) illustates
the procedure.

{\epsfxsize=200pt 

\makebox(300,110)[b]{\epsfbox{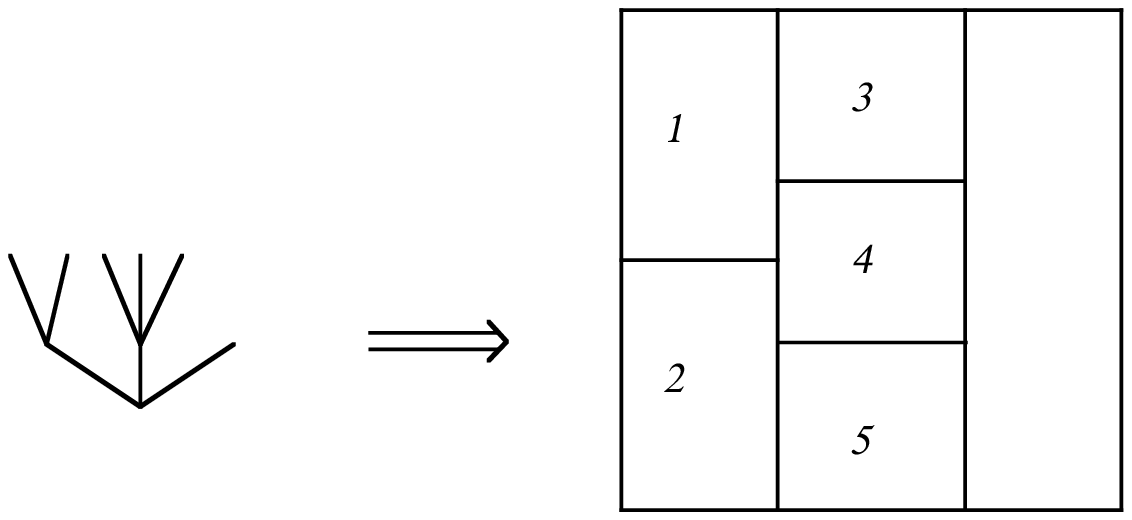}}}
\begin{figure}[h]\caption{}\label{Milgram}\end{figure}

Now for every morphism of trees $\sigma:T\rightarrow S$ we have to specify a homotopy $h_{\sigma}(t)$
$$h_{\sigma}(0) = \mu(a_S;a_{T_1},\ldots,a_{T_k})\ ; \  h_{\sigma}(1) = \pi(\sigma)a_T,$$ 
where $\mu$ is the multiplication in ${\cal C}^n$, $T_1,\ldots,T_k$ are fibers of $\sigma$,
and $\pi(\sigma)$ is the permutation which corresponds to $\sigma$ \cite{BEH}. 

This is more complicated combinatorially and we refer the reader to \cite{BEH}[Theorem 9.1] for the details
 in the case of $N(\M^n)$. The corresponding homotopies for ${\cal C}^n$ can be constructed analogously. 

For example, for the morphisms in $\Omega_2$ presented in Figure \ref{braiding}

  {\unitlength=0.71mm
\begin{picture}(100,30)(-15,0)
\put(15,5){\line(-1,1){5}}
\put(15,5){\line(1,1){5}}
\put(15,0){\line(0,1){5}}

\put(9,11){\line(1,1){5}}
\put(4,16){\line(0,1){5}}
\put(9,11){\line(-1,1){5}}

\put(21,11){\line(1,1){5}}
\put(26,16){\line(0,1){5}}
\put(21,11){\line(-1,1){5}}

\put(28,3){\vector(1,0){20}}

\put(49,0){\MT}

\put(93,3){\vector(-1,0){20}}

\put(105,5){\line(-1,1){5}}
\put(105,5){\line(1,1){5}}
\put(105,0){\line(0,1){5}}

\put(99,11){\line(1,1){5}}
\put(104,16){\line(0,1){5}}
\put(99,11){\line(-1,1){5}}

\put(111,11){\line(1,1){5}}
\put(106,16){\line(0,1){5}}
\put(111,11){\line(-1,1){5}}

\end{picture}}

\begin{figure}[h] \caption{}\label{braiding}\end{figure}

\noindent the corresponding homotopies are

{\epsfxsize=250pt 

\makebox(300,100)[b]{\epsfbox{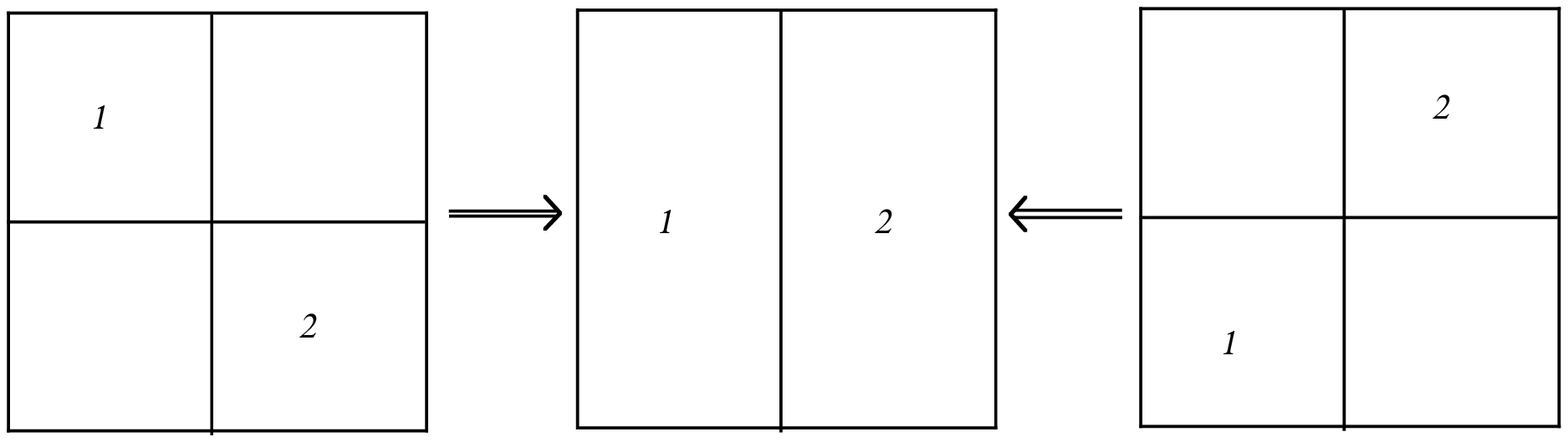}}} 
\begin{figure}[h]\caption{}\label{homotopy}
\end{figure}

The next piece of structure will be a homotopy
$$h_{\sigma_1,\sigma_2}:\Delta(2) \rightarrow {\cal C}^n_k$$
defined for every chain 
$$T\stackrel{\sigma_1}{\longrightarrow}S\stackrel{\sigma_2}{\longrightarrow}R$$
which agrees with the appropriate $h_{\sigma}$ on the faces of the $2$-simplex. This process continues.  
  
\

Theorem \ref{representable} shows that the combinatorial structure of an internal $n$-operad is actually the
quintessence of the algebraic structure which  $E_n$-operads were invented for. We can reformulate Theorem
\ref{representable} in the following way  
\begin{theorem} Every cofibrant $E_n$-operad $A$ is homotopically freely generated by its internal $n$-operad, in
the sense that there is a homotopy  equivalence of Kan simplicial sets
$${\mbox{\sc Oper}}_n(C) \simeq SOper^S(A,C).$$ 
An arbitrary $E_n$-operad is the representable object 
for the total left derived functor of ${\mbox{\sc Oper}}_n(-)$ on the homotopy category of
symmetric operads. \end{theorem}

\Remark According to this theorem every endomorphism of $E_n$-operads can be identified (up to higher homotopies) with an
internal operad in the little $n$-cube operad. This shows that the homotopy type of ${\mbox{\sc Oper}}_2({\cal
C}^2)$ is closely related to the Grothendieck-Teichm\"{u}ller group. It would be very interesting to study
these connections more closely.

\section{Coherence laws for $n$-fold loop spaces.}

It is  established in the previous section that the categorical $n$-operad 
 $\H^n$ provides us with all possible coherence laws. But this operad is very big  and in practice we only
need  some  suboperad or even quotient of a suboperad of $\H^1$. We will demonstrate this point 
with several  examples. 

\

\Example Let us begin from $n=1$. One can consider a suboperad ${\cal SH}^1$ of $\H^1$ which is
generated by the collection which contains one operation  in every strictly positive 
arity. So we do not consider  nullary operations. This suboperad is easy to describe:
the objects  are   chains of {\bf surjective} morphisms in $\Omega_1$ (see the description of $\H^n$ in  section
2.2).

Moreover, we want to have only one 
unary operation, so we take  quotion of this suboperad by the corresponding relation.
This means that we ask for the morphisms
\begin{equation}\label{TU}  (T\rightarrow U_n\rightarrow U_n) \longrightarrow (T\rightarrow U_n)\end{equation}
and
\begin{equation}\label{TT}(T\stackrel{id}{\rightarrow}T \rightarrow U_n)\longrightarrow (T \rightarrow U_n),\end{equation}
as well as their tensor products, to be identity morphisms in our factor operad. 

Then we have a poset operad which in arity $l$ is the poset of faces of the associahedron $K_l$ (see the
classical description of this poset in terms of planar trees \cite{Oper}). 

For example, for $l=3$ the picture of the corresponding poset is

{\unitlength=0.71mm
\begin{picture}(100,20)(-15,0)
\put(15,0){\line(-1,1){5}}
\put(15,0){\line(1,1){5}}
\put(10,6){\line(1,1){5}}
\put(20,6){\line(0,1){5}}
\put(10,6){\line(-1,1){5}}

\put(28,3){\vector(1,0){20}}

\put(60,0){\line(-1,1){5}}
\put(60,0){\line(1,1){5}}
\put(60,0){\line(0,1){5}}

\put(93,3){\vector(-1,0){20}}

\put(105,0){\line(-1,1){5}}
\put(105,0){\line(1,1){5}}
\put(110,6){\line(1,1){5}}
\put(100,6){\line(0,1){5}}
\put(110,6){\line(-1,1){5}}

\end{picture}}

\begin{figure}[h]\caption{}\end{figure}

The result for  $l=4$ is presented in Figure \ref{hexagon}.

\

\Remark The result of our identification is the contraction of the top face of the hexagon on Figure \ref{hexagon}. In general if
we do not require the tensor product of the morphisms (\ref{TU}),(\ref{TT}) to be the identity, we get the collection of
permutohedra
$P_{\star}$ and a quotient map  
$$q: P_{\star}\rightarrow K_{\star},$$
 which is of some importance as well \cite{Tonks,Loday}. The
collection of permutohedra does not have a structure of nonsymmetric operad so this quotient map is not a map of operads.
Nevertheless, we can define yet another  notion of {\it noncommutative nonsymmetric operad} by dropping the two  extreme 
conditions for associativity  in the definition of {\it Markl's pseudooperad} \cite{Oper}[Definition 1.16]. Every 
nonsymmetric operad is, of course, a noncommutative operad. The collection $P_{\star}$ is an example of a noncommutative 
operad and $q$ is an operadic map of non-commutative operads. 
 Recently W.Joyce used this notion of noncommutative operad implicitly in his work \cite{Joyce}.    

\

\Example Let us go to  higher dimensions.

 As with $n=1$  we do not consider nullary operations, so we consider a suboperad ${\cal SH}^n$ of $\H^n$ based on
morphisms of trees which are surjective in dimension $n$.  We also require  conditions  (\ref{TU}) and (\ref{TT}).  Yet, 
we want a little bit  more. 

The operad $\H^n$ has as generating objects all  trees including non-pruned trees. This corresponds to  possibility of 
having $l$-ary operations on  $n$-fold loop space which depend on `low dimensional units'. Generally speaking, this can be the
case in the nerve of an $n$-tuply monoidal $\omega$-category. In an algebra of the little $2$-cube operad,   for example, we  have
different binary multiplications which correspond to the following trees.

{\unitlength=0.71mm
\begin{picture}(100,20)(-15,0)
\put(15,0){\line(-1,1){5}}
\put(15,0){\line(1,1){5}}
\put(10,5){\line(1,1){5}}

\put(10,5){\line(-1,1){5}}

\put(60,0){\line(-1,1){5}}
\put(60,0){\line(1,1){5}}
\put(60,0){\line(0,1){5}}
\put(65,5){\line(1,1){5}}
\put(65,5){\line(-1,1){5}}

\put(105,0){\line(0,1){5}}

\put(105,5){\line(1,1){5}}
\put(105,5){\line(-1,1){5}}

\end{picture}}

\begin{center}{\sc Figure} \addtocounter{figure}{1}\thefigure: \end{center}

\noindent However, if  we are interested in the coherence laws for a unique binary multiplication on an $n$-fold loop space,
we want these low dimensional units to be strict. 

Formally we can do this as follows.
Observe, that the trees above has a common maximal pruned subtree. 
In general,
for a tree $T$ let $T^{(p)}$ be  its maximal
pruned subtree. We have a morphism $i^{(p)}:T^{(p)}\rightarrow T$ in
$\Omega_n$. It is not hard to see that $(-)^{(p)}$ can be extended to a functor (even an involutive comonad) on $\Omega_n$.
Therefore, it gives an  endofunctor on $\H^n$. Here we consider $\H^n$ as a category $$\coprod_{T} \H^n_T.$$   
Now we will factorise our  suboperad ${\cal SH}^n$ by the  equivalence relation generated by (\ref{TU}),(\ref{TT}) 
and  identification to the identity of those morphisms of trees $\sigma$ for which  $\sigma^{(p)}=
id$.

When $n=2$ the resulting category in the arity \MT has already appeared  in Figure \ref{braiding}.   

\begin{figure}[p]
{\epsfxsize=250pt 
\makebox(345,280)[t]{\epsfbox{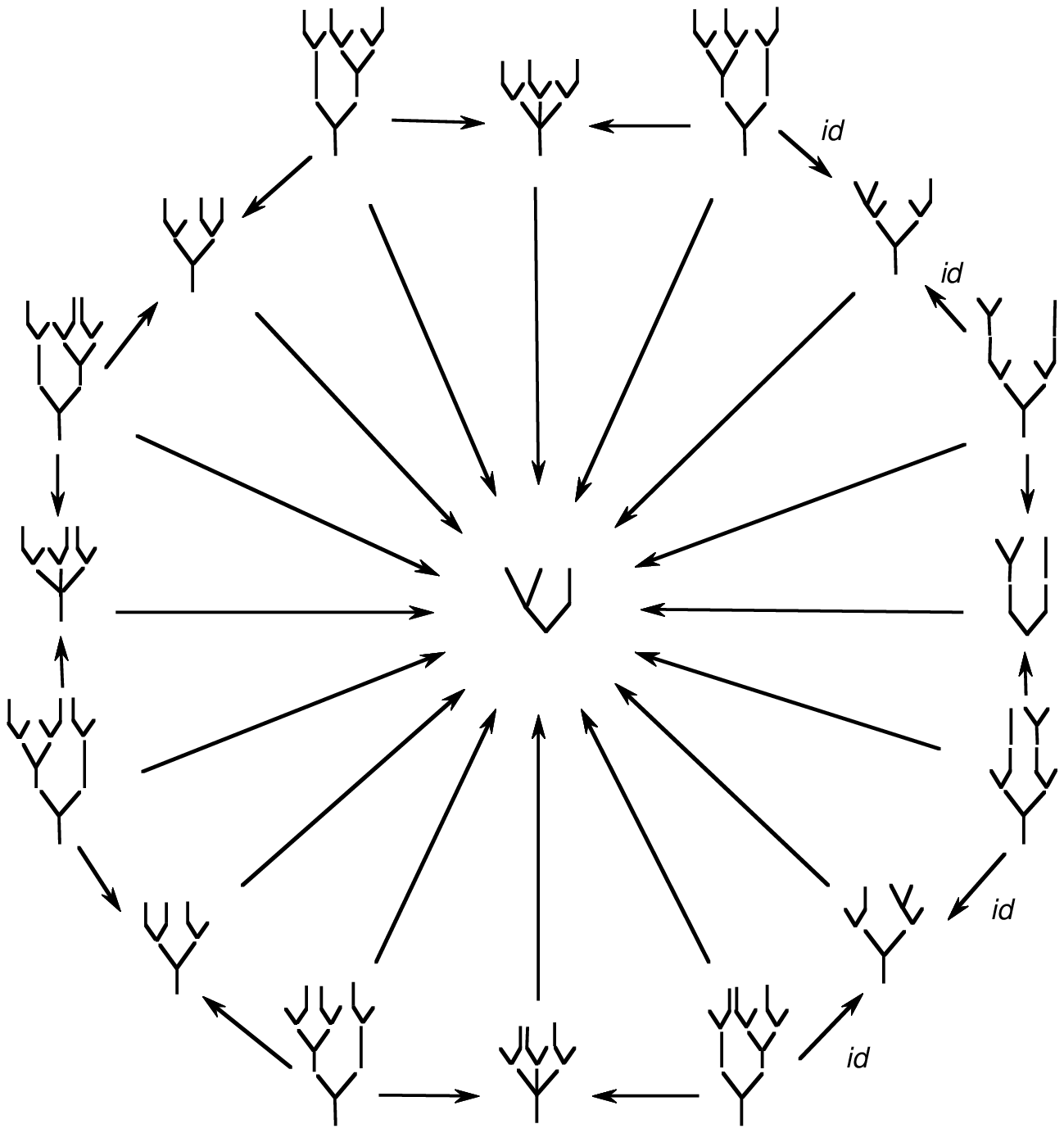}}}
\caption{}\label{octagon}\end{figure}

\noindent  Figure \ref{octagon} shows the corresponding category at the arity \MS. The reader can notice that this is
actually the poset of faces of one of the hexagons from the definition of the braided monoidal category
\cite{JS}. If we want the associativity of the multiplication to be strict we have to factorise 
further and we will obtain the first of the triangles from Figure \ref{triangle}. 
 
Using different degrees of  factorisation we can get many interesting polytopes which appear in the
literature like permutoassociahedra \cite{K}, resultohedra \cite{GKZ,KV} and their genearalisations.
We are going to study their properties in a future paper \cite{BS}.

\begin{figure}[p]
{\epsfxsize=220pt 
\makebox(340,180)[t]{\epsfbox{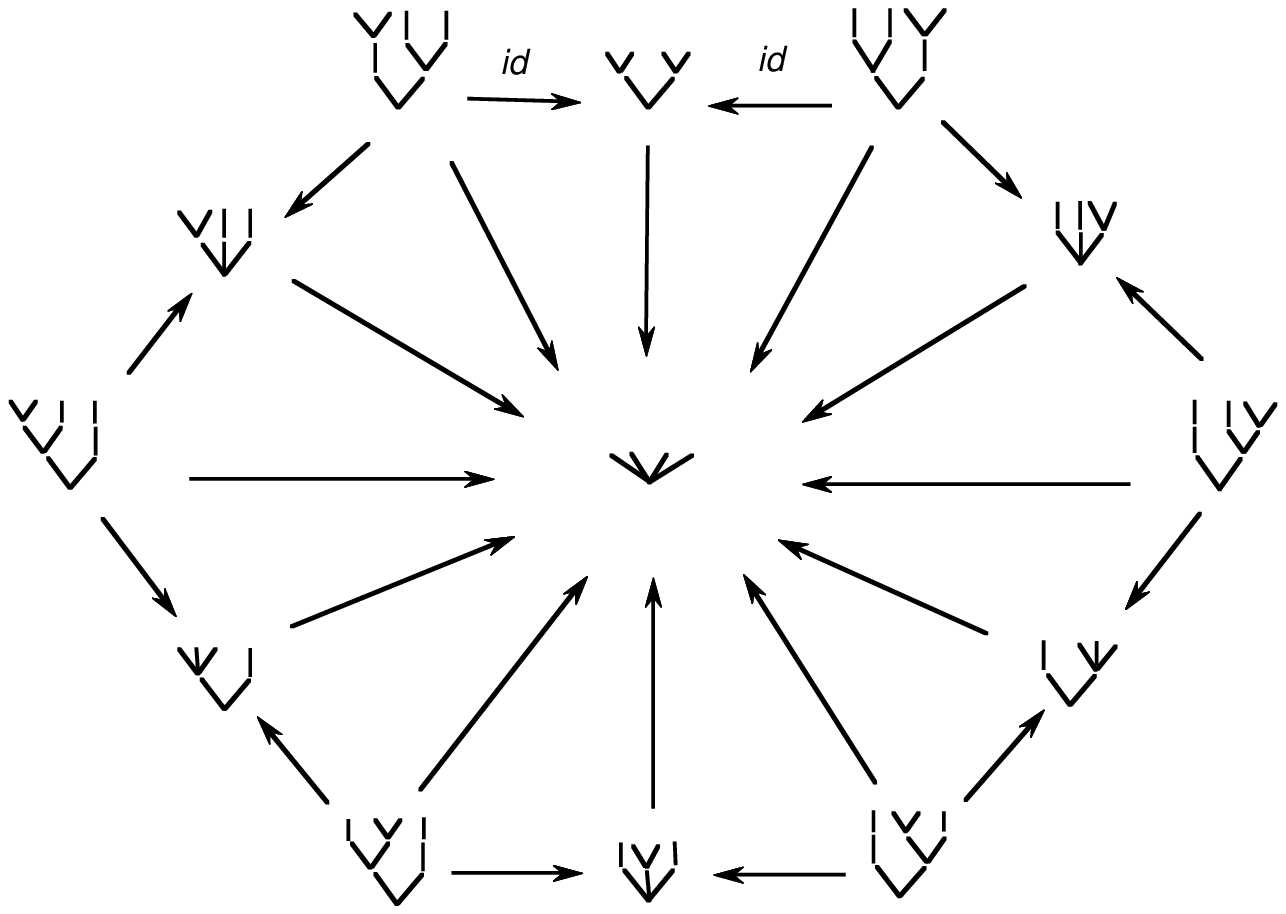}}}\caption{}\label{hexagon}\end{figure}

\newpage

\noindent {\bf Acknowledgements.} The first version of this paper was written during my visit to the Max-Plank Institut in
Bonn in  the automn of 2002 and I  gratefully acknowledge their support.   I also would like to thank C.Berger, A.Joyal,
A.Davydov, A.Makovecky, P.May, J.McClure,  J.Smith, and R. Street for many stimulating discussions during my work on this paper. 
 Finally, I acknowledge the financial support of the Scott Russell Johnson Memorial 
 Foundation and the
Macquarie University Research Commitee.

\end{document}